%% file: arxives_v2.tex
\documentclass{amsart}
\usepackage{amssymb,latexsym,amsmath,amscd,graphicx,color}

\theoremstyle{plain}
\newtheorem{thm}{Theorem}[section]
\newtheorem{prop}[thm]{Proposition}
\newtheorem{cor}[thm]{Corollary}
\newtheorem{lemma}[thm]{Lemma}

\theoremstyle{definition}
\newtheorem{defn}[thm]{Definition}

\newtheorem{rmks}[thm]{Remarks}
\newtheorem{notat}[thm]{Notation}

\newtheorem{exs}[thm]{Examples}

\newcommand{\lra}{\longrightarrow}

\newcommand{\NN}{\mathbb{N}}
\newcommand{\ZZ}{\mathbb{Z}}

\newcommand{\J}{\mathcal{J}}

\newcommand{\I}{\mathcal{I}}

\newcommand{\M}{\mathcal{M}}

\newcommand{\U}{\mathcal{U}}

\newcommand{\HH}{\mathcal{H}}
\newcommand{\X}{\mathcal{X}}

\newcommand{\LL}{\mathcal{L}}

\newcommand{\hgt}{\mbox{ht}\;}

\begin{document}

\title{Symmetric ladders and G-biliaison}

\author{Elisa Gorla}
\address{%
Institut f\"ur Mathematik\\
Universit\"at Z\"urich\\ 
Winterthurerstrasse 190\\
CH 8057 Z\"urich\\
Switzerland}
\email{elisa.gorla@math.uzh.ch}

\begin{abstract} 
We study the family of ideals generated by minors of
mixed size contained in a ladder of a symmetric matrix from the point
of view of liaison theory. We prove that they can be obtained from
ideals of linear forms by ascending G-biliaison. In particular, they
are glicci.
\end{abstract}

\thanks{The author was supported by the ``Forschungskredit der
  Universit\"at Z\"urich'' (grant no. 57104101) and by the Swiss National
  Science Foundation (grants no. 107887 and no. 123393).}

\subjclass{14M06, 13C40, 14M12}

\keywords{G-biliaison, Gorenstein liaison, 
minor, symmetric matrix, symmetric ladder, complete intersection, 
Gorenstein ideal, Cohen-Macaulay ideal.}

\maketitle

\section*{Introduction}

Ideals generated by minors have been studied extensively. They are
a central topic in commutative algebra, where they have been
investigated mainly using Gr\"obner bases and combinatorial
techniques (see among others~\cite{ea71}, \cite{kl79}, \cite{br81},
\cite{br88}, \cite{st90}, \cite{he92}). 
They are also relevant in algebraic geometry, since many classical
varieties such as the Veronese and the Segre variety are cut out by
minors. Degeneracy loci of morphisms between direct sums of line
bundles over projective space have a determinantal description, as do
the Schubert varieties. 

In this paper, we study ideals of minors in a symmetric
matrix from the point of view of liaison theory. In particular, we
consider ideals generated by minors of mixed size which are contained
in a symmetric ladder. Cogenerated ideals in a ladder of a symmetric
matrix belong to the family that we study. The family of cogenerated
ideals is a natural one to study from the combinatorial point of view
(see~\cite{de76} or~\cite{de82}). However, from the 
point of view of liaison theory it is more natural to study a larger
class of ideals, as they naturally arise during the linkage
process. We call them symmetric mixed ladder determinantal ideals.

In Section~1 we set the notation and define symmetric mixed ladder
determinantal ideals (Definition~\ref{ldi}). In Example~\ref{excog}
(3) we discuss why cogenerated ladder determinantal ideals of a
symmetric matrix are a special case of symmetric mixed ladder
determinantal ideals. In Proposition~\ref{primeCM} we show that
symmetric mixed ladder determinantal ideals are prime and
Cohen-Macaulay. In Proposition~\ref{ht} we express their height as the
cardinality of a suitable subladder.

In Section~2 we review the notion of G-biliaison, stating the
definition and main result in the algebraic language
(see Definition~\ref{bilid} and Theorem~\ref{bil}). In
Theorem~\ref{biliaison} we prove that symmetric mixed ladder determinantal
ideals can be obtained from ideals of linear forms by ascending
G-biliaison. In particular, they are glicci (Corollary~\ref{glicci}).

\section{Ideals of minors of a symmetric matrix}

Let $K$ be an algebraically closed field. 
Let $X=(x_{ij})$ be an $n\times n$ symmetric
matrix of indeterminates. In other words, the entries $x_{ij}$ with
$i\leq j$ are distinct indeterminates, and $x_{ij}=x_{ji}$ for
$i>j$. Let $K[X]=K[x_{ij}\mid 1\leq i\leq j\leq n]$ be the
polynomial ring associated to the matrix $X$. 
In this paper, we study ideals generated by the minors contained in a
ladder of a generic symmetric matrix from the point of view of liaison
theory. Throughout the paper, we only consider symmetric ladders. This
can be done without loss of generality, since the ideal generated by
the minors in a ladder of a symmetric matrix coincides with the ideal
generated by the minors in the smallest symmetric ladder containing it.

\begin{defn}\label{ladd}
A {\bf ladder} $\LL$ of $X$ is a subset of the set
$\X=\{(i,j)\in\NN^2 \mid 1\le i,j\le n\}$
with the following properties :
\begin{enumerate}
\item if $(i,j)\in\LL$ then $(j,i)\in \LL$ (i.e. $\LL$ is symmetric), and
\item if $i<h,j>k$ and $(i,j),(h,k)\in\LL$, then
$(i,k),(i,h),(h,j),(j,k)\in\LL$.
\end{enumerate}

We do not make any connectedness assumption on the ladder $\LL$.
For ease of notation, we also do not assume that $X$ is the 
smallest symmetric matrix containing $\LL$. 
Let $$\X^+=\{(i,j)\in\X\mid 1\leq i\leq j\leq n\}\;\;\;\mbox{and}\;\;\; 
\LL^+=\LL\cap\X^+.$$
Since $\LL$ is symmetric, $\LL^+$ determines $\LL$ and vice versa. We
will abuse terminology and call $\LL^+$ a ladder.
Observe that $\LL^+$ can be written as
$$\LL^+=\{(i,j)\in\X^+\mid i\leq c_l \mbox{ or } j\leq d_l \mbox{ for } 
l=1,\ldots,r\; \mbox{ and }$$ $$ i\geq a_l \mbox{ or } j\geq b_l 
\mbox{ for } l=1,\ldots,u \}$$ 
for some  integers $1\leq a_1<\ldots<a_u\leq n$, $n\geq b_1>\ldots>b_u\geq 1$, 
$1\leq c_1<\ldots<c_r\leq n$, and $n\geq d_1>\ldots>d_r\geq 1$, with
$a_l\leq b_l$ for $l=1,\ldots,u$ and $c_l\leq d_l$ for $l=1,\ldots,r$.

The points $(a_1,b_2),\ldots,(a_{u-1},b_u)$ are the {\bf lower outside
corners} of the ladder, $(a_1,b_1),\ldots,(a_u,b_u)$ are the
{\bf lower inside corners}, $(c_2,d_1),\ldots,(c_r,d_{r-1})$ the {\bf upper outside
corners}, and $(c_1,d_1),\ldots,(c_r,d_r)$ the {\bf upper inside
corners}. If $a_u\neq b_u$, then $(a_u,a_u)$ is a lower outside corner
and we set $b_{u+1}=a_u$. Similarly, if $c_r\neq d_r$ then $(d_r,d_r)$ is an
upper outside corner, and we set $c_{r+1}=d_r$. See also
Figure~\ref{fig1}. A ladder has at least one upper and one lower outside
corner. Moreover, $(a_1,b_1)=(c_1,d_1)$ is both an upper and a lower
inside corner.
\begin{figure}[htbp]
\begin{center}
\input{ladd1.pstex_t}
\caption{An example of ladder with tagged lower and upper corners.}
\label{fig1}
\end{center}
\end{figure}
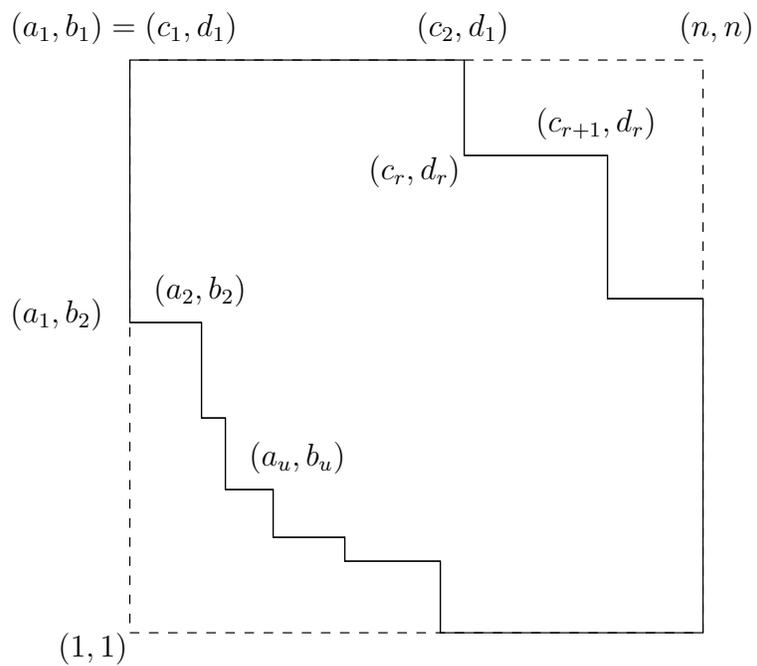

The {\bf upper border} of $\LL^+$ consists of the elements $(c,d)$ of $\LL^+$
such that either $c_l\leq c\leq c_{l+1}$ and $d=d_l$, or $c=c_l$ and
$d_l\leq d\leq d_{l-1}$ for some $l$. See Figure~\ref{fig2}.
\begin{figure}[htbp]
\begin{center}
\input{ladd2.pstex_t}
\caption{The upper border of the same ladder.}
\label{fig2}
\end{center}
\end{figure}
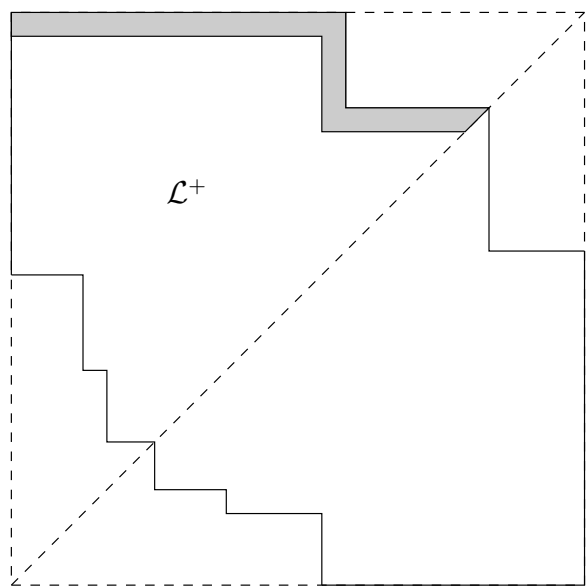
\end{defn}
All the corners belong to $\LL^+$. In fact, the ladder $\LL^+$
corresponds to its set of lower and upper outside (or equivalently
lower and upper inside) corners. The upper corners
of a ladder belong to its upper border.


Given a ladder $\LL$ we set $L=\{x_{ij}\in X\mid (i,j)\in\LL^+\}$,
and denote by $K[L]$ the polynomial ring $K[x_{ij}\mid x_{ij}\in
L]$. For $t$ a positive integer, and
$1\leq\alpha_1\leq\ldots\leq\alpha_t\leq n$,
$1\leq\beta_1\leq\ldots\leq\beta_t\leq n$ integers, we denote by  
$[\alpha_1,\ldots,\alpha_t; \beta_1,\ldots,\beta_t]$ the $t$-minor
$\det(x_{\alpha_i,\beta_j})$.
We let $I_t(L)$ denote the ideal generated by the set of the
$t$-minors of $X$ which involve only indeterminates of
$L$. In particular $I_t(X)$ is the ideal of $K[X]$ generated by 
the minors of $X$ of size $t\times t$.

In this article, we study the G-biliaison class of a large family of
ideals generated by minors in a ladder of a symmetric matrix.

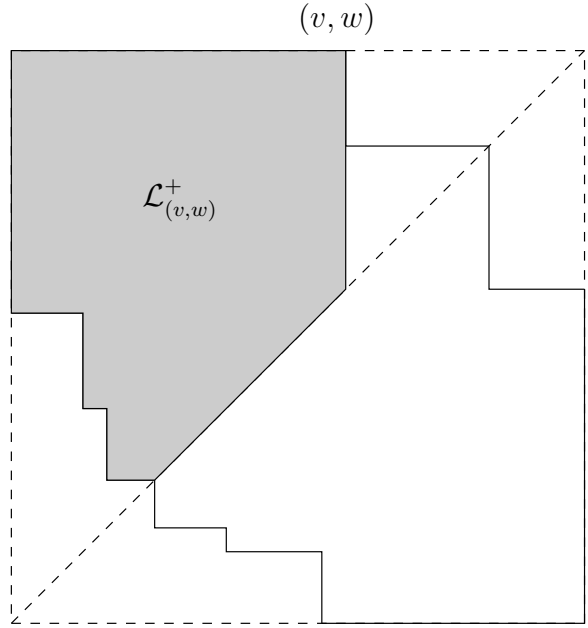
\begin{figure}[htbp]
\begin{center}
\input{ladd3b.pstex_t}
\caption{The ladder $\LL^+$ with a shaded subladder
  $\LL^+_{(v,w)}$.}
\label{fig3}
\end{center}
\end{figure}

\begin{notat}\label{decomp}
Let $\LL^+$ be a ladder. For $(c,d)\in\LL^+$ let 
$$\LL^+_{(c,d)}=\{(i,j)\in\LL^+\mid i\leq c,\; j\leq d\},
\;\;\;\;\; L_{(c,d)}=\{x_{ij}\in X\mid (i,j)\in\LL^+_{(c,d)}\}.$$
See also Figure~\ref{fig3}. 
Notice that $\LL^+_{(c,d)}$ is a ladder according to
Definition~\ref{ladd} and 
$$\LL^+=\bigcup_{(c,d)\in\U}\LL^+_{(c,d)}$$
where $\U$ denotes the set of upper outside corners of $\LL^+$.
\end{notat}

\begin{defn}\label{ldi}
Let $\{(v_1,w_1),\ldots,(v_s,w_s)\}$ be a subset of the
upper border of $\LL^+$ which contains all the upper outside
corners. We order them so that $1\leq v_1\leq\ldots\leq v_s\leq n$
and $n\geq w_1\geq\ldots\geq w_s\geq 1$.
Let $t=(t_1,\ldots,t_s)$ be a vector of positive
integers. Denote $L_{(v_k,w_k)}$ by $L_k$. The ideal 
$$I_t(L)=I_{t_1}(L_1)+\ldots+I_{t_s}(L_s)$$
is a {\bf symmetric mixed ladder determinantal ideal}. 
Denote $I_{(t,\ldots,t)}(L)$ by $I_t(L)$.
We call $(v_1,w_1),\ldots,(v_s,w_s)$ {\bf distinguished points} of
$\LL^+$.
\end{defn}
%

\begin{rmks}\label{conv}
\begin{enumerate}
\item Let  $\M\supseteq\LL$ be two ladders of $\X$, 
and let $M,L$ be the corresponding sets of indeterminates. 
We have isomorphisms of graded $K$-algebras 
$$K[L]/I_t(L)\cong K[M]/I_t(L)+(x_{ij}\;|\; x_{ij}\in
M\setminus L)\cong K[X]/I_{2t}(L)+(x_{ij}\;|\; x_{ij}\in
X\setminus L).$$
Here $I_t(L)$ is regarded as an ideal in $K[L],K[M],$ and $K[X]$
respectively. 
Then the height of the ideal $I_t(L)$ and the property of being prime, Cohen-Macaulay,
Gorenstein, Gorenstein in codimension $\leq c$ (see
Definition~\ref{gc}) do not depend on
whether we regard it as an ideal of $K[L],K[M],$ or $K[X]$. 
\item We can assume without loss of generality that for each
  $l=1,\ldots,s$ there exists a $k\in\{1,\ldots,u-1\}$ such that 
$$t_l\leq\min\{v_l-a_k+1, w_l-b_{k+1}+1\}$$ In fact, if 
$t_l>\min\{v_l-a_k+1, w_l-b_{k+1}+1\}$ for all $k$, then
$I_{t_l}(L_l)=0$. If that is the case, replace $L$ by
$M:=\cup_{i\neq l}L_i$, eliminate $(v_l,w_l)$ from the distinguished
points and remove the $l$-th entry of $t$ to get a new vector $m$. Then
we obtain a new ladder for which the assumption is satisfied and such
that $I_m(M)=I_t(L)$.
\item We can assume that 
$$w_k-w_{k-1}<t_k-t_{k-1}<v_k-v_{k-1},\;\;\;\mbox{for } k=2,\ldots,s.$$ 
In fact, if $v_k-v_{k-1}\leq t_k-t_{k-1}$, by successively developing a
$t_k$-minor of $L_k$ with respect to the first $v_k-v_{k-1}$ rows we
obtain an expression of the minor as a combination of minors of 
size $t_k-(v_k-v_{k-1})\geq t_{k-1}$ that involve only indeterminates
from $L_{k-1}$. Therefore $I_{t_k}(L_k)\supseteq I_{t_{k-1}}(L_{k-1})$.
Similarly, if $w_k-w_{k-1}\geq t_k-t_{k-1}$, by developing a $t_{k-1}$-minor
of $L_{k-1}$ with respect to the last $w_{k-1}-w_k$ columns
we obtain an expression of the minor as a combination of minors of 
size $t_{k-1}-(w_{k-1}-w_k)\geq t_k$ that involve only indeterminates
from $L_k$. Therefore $I_{t_{k-1}}(L_{k-1})\subseteq I_{t_k}(L_k)$.
In either case, we can remove a part of the ladder and reduce to the
study of a proper subladder that corresponds to the same symmetric
ladder determinantal ideal. 
\item We can always find $k\in\{1,\ldots,s\}$ such that $v_k>v_{k-1}$ and
$w_k>w_{k+1}$. In fact, the two inequalities are satisfied if and only
if $(v_k,w_k)$ is an upper outside corner. Notice that if we have
distinguished points $(v_k,w_k)$ and $(v_{k+1},w_{k+1})$ on the same
row or column, then one of the following holds: 
\begin{itemize}
\item either $v_k=v_{k+1}$ and $t_k>t_{k+1}$,
\item or $w_k=w_{k+1}$ and $t_{k+1}>t_k$.
\end{itemize}
In particular, we can find $k\in\{1,\ldots,s\}$ such that
$t_k\geq 2$, $v_k>v_{k-1}$ and $w_k>w_{k+1}$, unless $t_k=1$ for all
$k$.
\end{enumerate}
\end{rmks}

The following are examples of determinantal ideals of a symmetric
matrix which belong to the class of ideals that we study.

\begin{exs}\label{excog}
\begin{enumerate}
\item If $t=(t,\ldots,t)$ then $I_t(L)$ is the
ideal generated by the $t$-minors of $X$ that
involve only indeterminates from $L$. These ideals have been
studied in~\cite{co94a},  \cite{co94b}, and  \cite{co94c}.
\item If $\LL=\X$, then  according to Remarks~\ref{conv} we can assume
  that $w_l=n$ for all $l=1,\ldots,s$ and $v_s=n$. From
  Remark~\ref{conv} (3), we have $t_l>t_{l-1}$ and $v_l>v_{l-1}$ for
  all $l$. Then $I_t(L)$ is generated by the $t_1$-minors of the first $v_1$
  rows, the $t_2$-minors of the first $v_2$ rows, \dots, the
  $t_s$-minors of the whole matrix. This is a simple example of a
  cogenerated ideal.
\item The family of symmetric mixed ladder determinantal ideals contains
  the family of {\bf cogenerated ideals} in a ladder of a symmetric
  matrix, as defined
  in~\cite{co94a}. We follow the notation of~\cite{co94a}, and assume
  for ease of notation that $(1,n)$ is an inside corner of
  $\LL$ (i.e., that $X$ is the smallest matrix containing $L$). If
  $\alpha=\{\alpha_1,\ldots,\alpha_t\}$, then $I_{\alpha}(L)=I_{\tau}(L)$ where
  $\{(v_1,w_1),\ldots,(v_s,w_s)\}$ consists of the upper outside corners
  of $\LL$, together with the points of the upper border of $\LL$ which belongs to
  row $\alpha_l-1$, for all $l$ for which such an intersection point
  is unique (if for some $l$ the intersection of the row $\alpha_l-1$
  with the upper border of $\LL$ consists of more than one point, then
  $\LL$ has an upper outside corner on the row $\alpha_l-1$ and we do
  not add any extra point to the set). For each $k=1,\ldots,s$, we let
  $\tau_k=\min\{l\mid \alpha_l>v_k\}$.
\item Let $X$ be a matrix of size $m\times n$, $m\leq n$, whose
  entries are indeterminates. Assume that $X$ contains a square
  symmetric submatrix of indeterminates, and that all the other
  entries of $X$ are distinct indeterminates. In block
  notation 
$$X=\left(\begin{array}{cc} M & N \\ S & P 
\end{array}\right)$$
where $S$ is a symmetric matrix of indeterminates and $M,N,P$ are
generic matrices of indeterminates. Let $t\in\ZZ_+$. Then $I_t(X)$ is
a symmetric ladder determinantal ideal generated by the minors of size
$t\times t$ contained in a symmetric ladder of 
$$\left(\begin{array}{ccc} Y & M & N \\ 
M^t & S & P \\ N^t & P^t & Z\end{array}\right)$$ where $Y,Z$ are
symmetric matrices of indeterminates, and $M^t$ denotes the transpose
of $M$. This was observed by Conca in~\cite{co94a}.
\end{enumerate}
\end{exs}

In this section we establish some properties of symmetric mixed ladder
determinantal ideals. It is known (\cite{co94a}) that cogenerated
ideals are prime and Cohen-Macaulay. In the sequel we show that the
result of Conca easily extends to symmetric mixed ladder determinantal
ideals. We exploit a well known localization technique
(see~\cite{br93}, Lemma~7.3.3). The same argument was used to prove
Lemma~1.19 in~\cite{go07b}. For completeness, we state it for the case of a ladder
of a symmetric matrix and we outline the proof. We use
the notation of Definitions~\ref{ladd} and~\ref{ldi}. From
Remark~\ref{conv} (4) we know that we can always find
$k\in\{1,\ldots,s\}$ such that $t_k\geq 2$, $v_k>v_{k-1}$ and
$w_k>w_{k+1}$, unless $t=(1,\ldots,1)$.

\begin{lemma}\label{local}
Let $\LL$ be a ladder of a symmetric matrix $X$ of
indeterminates. $\LL$ has a set of distinguished points
$\{(v_1,w_1),\ldots,(v_s,w_s)\}\in\LL^+$ and
$t=(t_1,\ldots,t_s)\in\ZZ_+^s$. 
Let $I_t(L)$ be the corresponding symmetric mixed ladder determinantal ideal.
Let $k\in\{1,\ldots,s\}$ such that $t_k\geq 2$, $v_k>v_{k-1}$ and
$w_k>w_{k+1}$.

Let $t'=(t_1,\ldots,t_{k-1},t_k-1,t_{k+1},\ldots,t_s)$ and 
let $\LL'$ be the ladder obtained from $\LL$ by removing the entries
$(v_{k-1}+1,w_k),\ldots,(v_k-1,w_k),(v_k,w_k),(v_k,w_k-1)\ldots, 
(v_k,w_{k+1}+1)$ and the symmetric ones. Let
$$(v_1,w_1),\dots,(v_{k-1},v_{k-1}),(v_k-1,w_k-1),
(v_{k+1},w_{k+1}),\dots,(v_s,w_s)$$ be the distinguished points of $\LL'$.

Then there is an isomorphism between $K[L]/I_t(L)[x_{v_k,v_k}^{-1}]$ and 
$$K[L']/I_{t'}(L')[x_{v_{k-1}+1,w_k},\ldots,x_{v_k-1,w_k},
x_{v_k,w_k}^{\pm 1},x_{v_k,w_k-1},\ldots, x_{v_k,w_{k+1}+1}].$$
\end{lemma}

\begin{proof}
Under the assumption of the lemma, $\LL'$ is a ladder and $I_{t'}(L')$
is a symmetric mixed ladder determinantal ideal. 
Let $$A=K[L][x_{v_k,w_k}^{-1}]$$ and
$$B=K[L'][x_{v_{k-1}+1,w_k},\ldots,x_{v_k-1,w_k},
x_{v_k,w_k}^{\pm 1},x_{v_k,w_k-1},\ldots, x_{v_k,w_{k+1}+1}].$$ 
Define a $K$-algebra homomorphism
$$\begin{array}{lcl}\varphi:A & \lra & B \\
x_{i,j} & \longmapsto & \left\{\begin{array}{ll} 
x_{i,j}+x_{i,w_k}x_{v_k,j}x_{v_k,w_k}^{-1} & \mbox{if $i\neq v_k,
  j\neq w_k$ and $(i,j)\in\LL_{(v_k,w_k)}$} \\
x_{i,j} & \mbox{otherwise.}
\end{array}\right.\end{array}$$
The inverse of $\varphi$ is 
$$\begin{array}{lcl}\psi:B & \lra & A \\
x_{i,j} & \longmapsto &\left\{\begin{array}{ll} 
x_{i,j}-x_{i,w_k}x_{v_k,j}x_{v_k,w_k}^{-1} & \mbox{if $i\neq v_k,
  j\neq w_k$ and $(i,j)\in\LL'_{(v_k,w_k)}$} \\
x_{i,j} & \mbox{otherwise.}
\end{array}\right.\end{array}$$
It is easy to check that $\varphi$ and $\psi$ are inverse to each
other. Since
$$\varphi(I_{t_k}(L_{(v_k,w_k)})A)=I_{t_k-1}(L'_{(v_k-1,w_k-1)})B$$ we
have $$\varphi(I_t(L)A)=I_{t'}(L')B\;\;\; \mbox{hence}\;\;\;
A/I_t(L)A\cong B/I_{t'}(L')B.$$
\end{proof}

Using Lemma~\ref{local} we can establish some
properties of symmetric mixed ladder determinantal ideals. 

\begin{prop}\label{primeCM}
Symmetric mixed ladder determinantal ideals are prime and Cohen-Macaulay.
\end{prop}

\begin{proof}
Let $I_t(L)$ be the symmetric mixed ladder determinantal ideal associated to the
ladder $\LL$ with distinguished points $(v_1,w_1),\ldots,(v_s,w_s)$
and $t=(t_1,\ldots,t_s)$. Let $t_{\mbox{\tiny
    max}}=\max\{t_1,...,t_s\}$. If $t_{\mbox{\tiny
    max}}=1$ then $I_t(L)$ is generated by indeterminates, hence it is
prime and Cohen-Macaulay. Therefore assume that $t_{\mbox{\tiny
    max}}\geq 2$ and let ${\widehat \LL}$ be the ladder with the same
lower outside corners as $\LL$, and upper outside corners 
$(v_k+t_{\mbox{\tiny max}}-t_k,w_k+t_{\mbox{\tiny max}}-t_k)$ for
$k=1,\ldots,s$. Notice that the corners are distinct, and the
inequalities of Definition~\ref{ladd} are satisfied by
Remark~\ref{conv} (3). In other words, for each $k=2,\ldots,s$ we have 
$$w_k+t_{\mbox{\tiny max}}-t_k<w_{k-1}+t_{\mbox{\tiny max}}-t_{k-1}\;\;\;
\mbox{and}\;\;\; 
v_k+t_{\mbox{\tiny max}}-t_k>v_{k-1}+t_{\mbox{\tiny max}}-t_{k-1}.$$
Let ${\widehat L}=\{x_{ij}\in X\mid (i,j)\in
{\widehat\LL},\; i\leq j\}$ and let $M={\widehat L}\setminus L$.
Denote by $I_{t_{\mbox{\tiny max}}}({\widehat L})$ the ideal generated
by the minors of size $t_{\mbox{\tiny max}}$ which involve only
indeterminates in ${\widehat L}$. By Lemma~\ref{local} there exists 
a subset $\{z_1,...,z_m\}$ of $M$
such that $$K[{\widehat L}]/I_{t_{\mbox{\tiny max}}}({\widehat
  L})[z_1^{-1},\ldots, z_m^{-1}]\cong 
K[L]/I_t(L)[M][z_1^{-1},\ldots, z_m^{-1}].$$
The ring $K[{\widehat L}]/I_{t_{\mbox{\tiny max}}}({\widehat L})$ is a
Cohen-Macaulay domain by Theorem~1.13 in~\cite{co94a}. Therefore
$K[L]/I_t(L)[M][z_1^{-1},\ldots, z_m^{-1}]$ is a Cohen-Macaulay
domain. Since $M$ is a set of indeterminates over the ring $K[L]/I_t(L)$ and
$z_1,....,z_m\in M$, then $K[L]/I_t(L)[M]$ is a
Cohen-Macaulay domain. Hence  $I_t(L)$ is prime and Cohen-Macaulay.
\end{proof}

\begin{figure}[htbp]
\begin{center}
\input{ladd4b.pstex_t}
\caption{An example of $\LL^+$ with three distinguished points and
  $t=(3,6,4)$. The corresponding $\HH^+$ is shaded.}
\label{fig4}
\end{center}
\end{figure}
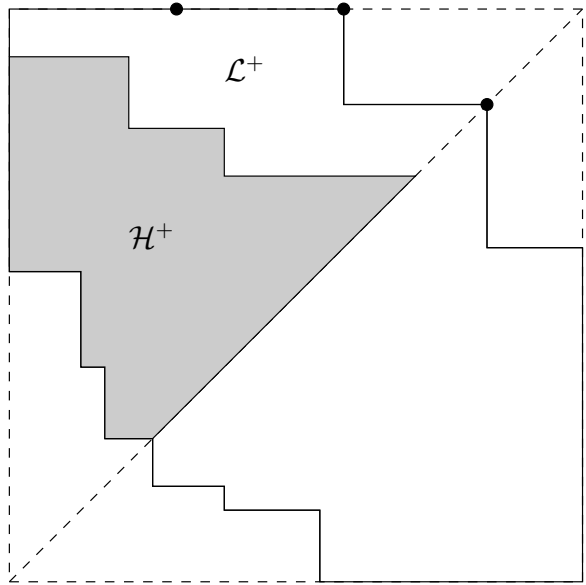

A standard argument allows us to compute the height of
symmetric mixed ladder determinantal ideals. These heights have been
computed by Conca in~\cite{co94a} for the family of cogenerated
ideals. The arguments in~\cite{co94a} are of a more combinatorial
nature, and the height is expressed as a sum of lengths of
maximal chains in some subladders. Our formula for the height is very
simple. The proof is independent of the results of Conca, and it
essentially follows from Lemma~\ref{local}.
We use the same notation as in Definitions~\ref{ladd}
and~\ref{ldi}, and Lemma~\ref{local}. An example is given in
Figure~\ref{fig4}.

\begin{prop}\label{ht}
Let $\LL$ be a ladder with distinguished points
$(v_1,w_1),\ldots,(v_s,w_s)$ and let  
$$\HH^+=\{(i,j)\in\LL^+\mid 
i\leq v_{k-1}-t_{k-1}+1 \mbox{ or } j\leq w_k-t_k+1 \mbox{ for }
k=2,\ldots,s,$$ $$j\leq w_1-t_1+1,\; i\leq v_s-t_s+1\}.$$ 
Let $\HH=\HH^+\cup\{(j,i)\mid (i,j)\in\HH^+\}.$
Then $\HH$ is a symmetric ladder and 
$$\hgt I_t(L)=|\HH^+|.$$
\end{prop}

\begin{proof}
Observe that by Remark~\ref{conv}
(3) $$v_k-t_k+1>v_{k-1}-t_{k-1}+1,\;\;\;\mbox{and}\;\;\; 
w_k-t_k+1<w_{k-1}-t_{k-1}+1.$$ 
Therefore $\HH$ is a ladder with upper outside corners 
$\{(v_k-t_k+1,w_k-t_k+1)\;|\;k=1,\ldots,s\}$ and the same lower
outside corners as $\LL$. Let $H=\{x_{i,j}\;|\; (i,j)\in\HH^+\}$.
We argue by induction on $\tau=t_1+\ldots+t_s\geq s$. If $\tau=s$,
then $t_1=\ldots=t_s=1$, and $\LL=\HH$. Hence
$$I_1(L)=(x_{ij}\mid x_{ij}\in L)=I_1(H)=|\HH^+|.$$ 
Assume now that the thesis holds for $\tau-1\geq s$ and prove it for
$\tau$. Since $\tau>s$, by Remark~\ref{conv} (4) there exists
$k\in\{1,\ldots,s\}$ such that $t_k\geq 2$, $v_k>v_{k-1}$ and
$w_k>w_{k+1}$. By Lemma~\ref{local} we have an
isomorphism between $K[L]/I_t(L)[x_{v_k,w_k}^{-1}]$ and 
$$K[L']/I_{t'}(L')[x_{v_{k-1},w_k},\ldots,x_{v_k-1,w_k},x_{v_k,w_k}^{\pm
1},x_{v_k,w_k-1},\ldots,
x_{v_k,w_{k+1}+1}].$$
Since $x_{v_k,w_k}$ does not divide zero modulo $I_{t'}(L')$ and $I_{t}(L)$, we have
$$\hgt I_{t}(L)=\hgt I_{t'}(L'). $$  
The thesis follows by the induction hypothesis, observing that
the same ladder $\HH$ computes the height of both $I_{t'}(L')$ and
$I_{t}(L)$. 
\end{proof}

\section{G-biliaison of symmetric mixed ladder determinantal ideals}

In this section we study symmetric mixed ladder determinantal ideals from
the point of view of liaison theory. We prove that they belong to the
G-biliaison class of a complete intersection. In particular, they are
glicci. This is yet another family of ideals of minors for which one
can perform a descending G-biliaison to an ideal in the same family,
in such a way that one eventually reaches an ideal generated by linear
forms. Other families of ideals that were treated with an analogous
technique are ideals generated by maximal minors of a matrix with polynomial
entries \cite{ha07}, minors of a symmetric matrix with
polynomial entries \cite{go07a}, minors of a matrix with polynomial
entries \cite{go07u}, minors of mixed size in a ladder of a generic
matrix \cite{go07b}, and pfaffians of mixed size in a ladder of a
generic skew-symmetric matrix \cite{de08u}.

In \cite{ha92}, \cite{ha94}, \cite{ha07} Hartshorne developed the
theory of generalized divisors, which is a useful language for
the study of Gorenstein liaison via the study of G-biliaison
classes. In~\cite{ha94} it was shown that even CI-liaison and
CI-biliaison generate the same equivalence classes. In~\cite{kl01}
Kleppe, Migliore, Mir\'o-Roig, Nagel and Peterson
proved that a G-biliaison on an arithmetically Cohen-Macaulay,
$G_1$ scheme can be realized via two G-links. The result was
generalized in~\cite{ha07} by Hartshorne to G-biliaison on an
arithmetically Cohen-Macaulay, $G_0$ scheme.

In Proposition~\ref{primeCM} we saw that symmetric mixed ladder
determinantal ideals are prime, hence they define reduced and
irreducible, projective algebraic varieties. Since we wish to work in
the algebraic setting, we state the definition of G-biliaison and the
main theorem connecting G-biliaison and G-liaison in the language of
ideals.

\begin{defn}\label{gc}
Let $R=K[L]$ and let $J\subseteq R$ be a homogeneous, saturated ideal. 
We say that $J$ is {\bf Gorenstein in codimension $\leq$ c} if the
localization $(R/J)_P$ is a Gorenstein ring for any prime ideal $P$ of
$R/J$ of height smaller than or equal to $c$. We often say that $J$ is
$G_c$. We call {\bf generically Gorenstein}, or $G_0$, an ideal $J$ which is
Gorenstein in codimension 0.
\end{defn}

\begin{defn}(\cite{ha07}, Sect. 3)\label{bilid}
Let $R=K[X]$ and let $I_1$ and $I_2$ be homogeneous ideals in
$R$ of pure height $c$. We say that $I_1$ is obtained by an {\bf
elementary G-biliaison} of height $h$ from $I_2$ if there exists a
Cohen-Macaulay, generically Gorenstein ideal $J$ in $R$ of height $c-1$
such that $J\subseteq I_1\cap I_2$ and $I_1/J\cong [I_2/J](-h)$ as
$R/J$-modules. If $h>0$ we speak about {\bf ascending} elementary
G-biliaison.
\end{defn}


The following theorem gives a connection between G-biliaison and
G-liaison.

\begin{thm}\label{bil}[Kleppe, Migliore, Mir\`o-Roig, Nagel,
  Peterson~\cite{kl01}; Hartshorne~\cite{ha07}]
Let $I_1$ be obtained by an elementary G-biliaison from $I_2$. Then
$I_2$ is G-linked to $I_1$ in two steps.
\end{thm}

We now show that symmetric mixed ladder determinantal ideals belong to the
G-biliaison class of a complete intersection. The idea of the proof is as
follows: starting from a symmetric mixed ladder determinantal ideal
$I$, we construct two symmetric mixed ladder determinantal ideals $I'$
and $J$ such that $J$ is contained in $I\cap I'$ and $\hgt I=\hgt
I'=\hgt J+1$. We show that $I$ can be obtained from $I'$ by an
elementary G-biliaison of height 1 on $J$.

\begin{thm}\label{biliaison}
Any symmetric mixed ladder determinantal ideal can be obtained from an ideal
generated by linear forms by a finite sequence of ascending elementary
G-biliaisons.
\end{thm}

\begin{proof}
Let $I_t(L)$ be a symmetric mixed ladder determinantal ideal associated to
a ladder $\LL^+$ with distinguished points
$(v_1,w_1),\ldots,(v_s,w_s)$. Let $L_k=L_{(v_k,w_k)}$, then
$$I_t(L)=I_{t_1}(L_1)+\cdots+I_{t_s}(L_s)\subseteq K[L].$$
As discussed in Remark~\ref{conv} (1) we will not distinguish between
symmetric mixed ladder determinantal ideals and their extensions. Therefore,
all ideals will be in $R=K[L]$.
If $t_1=\ldots=t_s=1$ then $I_t(L)$ is generated by linear
forms. Hence let $t_k=\max\{t_1,\ldots,t_s\}\geq 2$. From
Remark~\ref{conv} (3) we have that 
$w_{k+1}-w_k<0<v_k-v_{k-1}.$ In particular $(v_k,w_k)$ is an upper
outside corner.

Let $\LL'^+$ be the ladder with distinguished points
$$(v_1,w_1),\ldots,(v_{k-1},w_{k-1}),(v_k-1,w_k-1),(v_{k+1},w_{k+1}),\ldots,(v_s,w_s).$$
Observe that $\LL'^+$ is obtained from $\LL^+$ by removing the
entries $$(v_{k-1}+1,w_k),\ldots,(v_k-1,w_k),(v_k,w_k),(v_k,w_k+1),
\ldots,(v_k,w_{k+1}-1).$$
Let $t'=(t_1,\ldots,t_{k-1},t_k-1,t_{k+1},\ldots,t_s)$, and let
$I_{t'}(L')$ be the associated symmetric mixed ladder determinantal
ideal. 
It is easy to check that $\LL'^+$ and $t'$ satisfy
the inequalities of Definition~\ref{ldi} and of Remarks~\ref{conv}.
By Proposition~\ref{ht} $$\hgt I_t(L)=\hgt I_{t'}(L')=|\HH^+|$$ 
where $$\HH^+=\{(i,j)\in\LL^+\mid 
i\leq v_{k-1}-t_{k-1}+1 \mbox{ or } j\leq w_k-t_k+1 \mbox{ for }
k=2,\ldots,s,$$ $$j\leq w_1-t_1+1,\; i\leq v_s-t_s+1\}.$$

Let $\J^+$ be the ladder obtained from $\LL^+$ by removing 
$(v_k,w_k)$, and
let $$(v_1,w_1),\ldots,(v_{k-1},w_{k-1}),(v_k-1,w_k),(v_k,w_k-1),(v_{k+1},w_{k+1}),\ldots,(v_s,w_s)$$
be its distinguished points (see Figure~\ref{fig5}).
\begin{figure}[htbp]
\begin{center}
\input{ladd5.pstex_t}
\caption{An example of $\LL^+$ with $\LL'^+$ and $\J^+$. The
  distinguished point $(v_k,w_k)$ is marked. $\LL'^+$ is colored in a
  darker shade and the entries which belong to $\J^+$ but not to
  $\LL'^+$ are colored in a lighter shade.}
\label{fig5}
\end{center}
\end{figure}
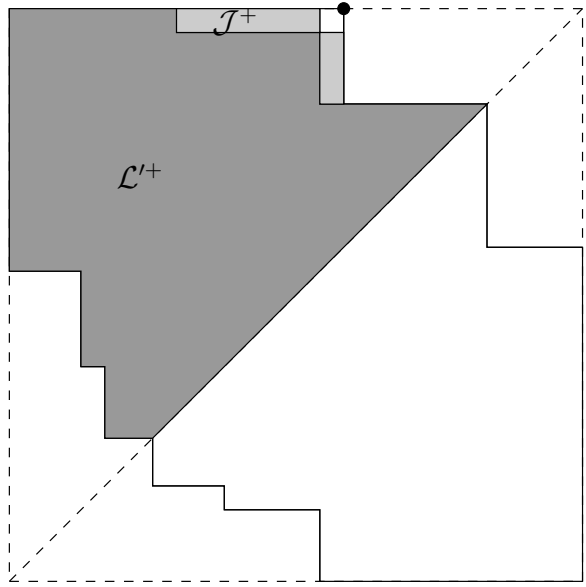
Let $u=(t_1,\dots,t_{k-1},t_k,t_k,t_{k+1},\dots,t_s)$. Then 
$$I_u(J)=I_{t_1}(L_1)+\ldots+I_{t_{k-1}}(L_{k-1})+I_{t_k}(J_{(v_k-1,w_k)})+
I_{t_k}(J_{(v_k,w_k-1)})+$$
$$+I_{t_{k+1}}(L_{k+1})+\ldots+I_{t_s}(L_s).$$
In other words, $I_u(J)$ is the ideal generated by the minors of
$I_t(L)$ that do not involve the indeterminate $x_{v_k,w_k}$.
We claim that $I_u(J)\subseteq I_t(L)\cap I_{t'}(L').$ It is clear
that $I_u(J)\subseteq I_t(L)$. The inclusion 
$I_u(J)\subseteq I_{t'}(L')$ follows from
$I_{t_k}(L_{(v_k-1,w_k)})+I_{t_k}(L_{(v_k,w_k-1)})\subset
I_{t_k-1}(L'_{(v_k-1,w_k-1)}).$

Let $$\I^+=\HH^+\setminus\{(v_k-t_k+1,w_k-t_k+1)\}.$$
By Proposition~\ref{ht} $$\hgt I_u(J)=|\I^+|=\hgt I_t(L)-1.$$ 
The ideal $I_u(J)$ is prime and Cohen-Macaulay by
Proposition~\ref{primeCM}. In particular it is generically
Gorenstein. 

We claim that $I_t(L)$ is obtained from $I_{t'}(L')$ by an
elementary G-biliaison of height $1$ on $I_u(J)$. This is equivalent
to showing that
\begin{equation}\label{isom}
I_t(L)/I_u(J)\cong [I_{t'}(L')/I_u(J)](-1)
\end{equation}
as $R/I_u(J)$-modules. Denote by
$[\alpha_1,\ldots,\alpha_t;\beta_1,\ldots,\beta_t]$
the $t\times t$-minor of $X$ which involves rows
$\alpha_1,\ldots,\alpha_t$ and columns
$\beta_1,\ldots,\beta_t$. We claim that multiplication by
$$f=\frac{[v_k-t_k+1,\ldots,v_k-1;
  w_k-t_k+1,\ldots,w_k-1]}{[v_k-t_k+1,\ldots,v_k-1,v_k;
  w_k-t_k+1,\ldots,w_k-1,w_k]}$$ yields an isomorphism
between $I_t(L)/I_u(J)$ and $[I_{t'}(L')/I_u(J)](-1)$.

Notice in fact that the ideal $I_t(L)/I_u(J)$
is generated by the minors of size
$t_k\times t_k$ of $L_k$ which involve both row $v_k$ and column
$w_k$, while the ideal $I_{t'}(L')/I_u(J)$
is generated by the minors of size $(t_k-1)\times (t_k-1)$ of $L'_k$.
For any minor
$[\alpha_1,\ldots,\alpha_{t_k-1},v_k;\beta_1,\ldots,\beta_{t_k-1},w_k]\in
I_{t_k}(L_k)$ which involves both row $v_k$ and column
$w_k$, consider the minor
$[\alpha_1,\ldots,\alpha_{t_k-1};\beta_1,\ldots,\beta_{t_k-1}]\in
I_{t_k-1}(L'_k)$. 
By \cite{go07a}, Lemma~2.6 
$$[\alpha_1,\ldots,\alpha_{t_k-1};\beta_1,\ldots,\beta_{t_k-1}]\cdot
[v_k-t_k+1,\ldots,v_k-1,v_k;w_k-t_k+1,\ldots,w_k-1,w_k]=$$
$$[v_k-t_k+1,\ldots,v_k-1;w_k-t_k+1,\ldots,w_k-1]\cdot
[\alpha_1,\ldots,\alpha_{t_k-1},v_k;\beta_1,\ldots,\beta_{t_k-1},w_k]$$
modulo $I_u(J)$.
Therefore the
ideals $$[v_k-t_k+1,\ldots,v_k-1;w_k-t_k+1,\ldots,w_k-1]\cdot I_t(L)+I_u(J)$$
and 
$$[v_k-t_k+1,\ldots,v_k;w_k-t_k+1,\ldots,w_k]\cdot I_{t'}(L')+I_u(J)$$ are
equal, hence they are equal modulo $I_u(J)$.
Therefore isomorphism~(\ref{isom}) holds, and $I_t(L)$ and
$I_{t'}(L')$ are G-bilinked on $I_u(J)$. Repeating this procedure, one
eventually reaches the ideal generated by the entries of the ladder
$\HH$ defined in Proposition~\ref{ht}. Clearly $$I_1(H)=(x_{ij}\;|\;
(i,j)\in\HH)$$ is a complete intersection.
\end{proof}

The following is a straightforward consequence of
Theorem~\ref{biliaison}, according to Theorem~\ref{bil}.

\begin{cor}\label{glicci}
Every symmetric mixed ladder determinantal ideal $I_t(L)$ can be G-linked in
$2(t_1+\ldots+t_s-s)$ steps to a complete intersection of linear forms
of the same height. Hence symmetric mixed ladder determinantal ideals are
glicci.
\end{cor}

\end{document}

%% file: ladd1.pstex_t
\begin{picture}(0,0)%
\includegraphics{ladd1.pstex}%
\end{picture}%
\setlength{\unitlength}{3947sp}%
\begingroup\makeatletter\ifx\SetFigFont\undefined%
\gdef\SetFigFont#1#2#3#4#5{%
  \reset@font\fontsize{#1}{#2pt}%
  \fontfamily{#3}\fontseries{#4}\fontshape{#5}%
  \selectfont}%
\fi\endgroup%
\begin{picture}(4769,4114)(451,-4169)
\put(3001,-211){\makebox(0,0)[lb]{\smash{{\SetFigFont{12}{14.4}{\familydefault}{\mddefault}{\updefault}{\color[rgb]{0,0,0}$(c_2,d_1)$}%
}}}}
\put(451,-2011){\makebox(0,0)[lb]{\smash{{\SetFigFont{12}{14.4}{\rmdefault}{\mddefault}{\updefault}{\color[rgb]{0,0,0}$(a_1,b_2)$}%
}}}}
\put(2701,-1111){\makebox(0,0)[lb]{\smash{{\SetFigFont{12}{14.4}{\rmdefault}{\mddefault}{\updefault}{\color[rgb]{0,0,0}$(c_r,d_r)$}%
}}}}
\put(451,-211){\makebox(0,0)[lb]{\smash{{\SetFigFont{12}{14.4}{\rmdefault}{\mddefault}{\updefault}{\color[rgb]{0,0,0}$(a_1,b_1)=(c_1,d_1)$}%
}}}}
\put(1351,-1861){\makebox(0,0)[lb]{\smash{{\SetFigFont{12}{14.4}{\rmdefault}{\mddefault}{\updefault}{\color[rgb]{0,0,0}$(a_2,b_2)$}%
}}}}
\put(1951,-2911){\makebox(0,0)[lb]{\smash{{\SetFigFont{12}{14.4}{\rmdefault}{\mddefault}{\updefault}{\color[rgb]{0,0,0}$(a_u,b_u)$}%
}}}}
\put(3751,-811){\makebox(0,0)[lb]{\smash{{\SetFigFont{12}{14.4}{\rmdefault}{\mddefault}{\updefault}{\color[rgb]{0,0,0}$(c_{r+1},d_r)$}%
}}}}
\put(4651,-211){\makebox(0,0)[lb]{\smash{{\SetFigFont{12}{14.4}{\rmdefault}{\mddefault}{\updefault}{\color[rgb]{0,0,0}$(n,n)$}%
}}}}
\put(751,-4111){\makebox(0,0)[lb]{\smash{{\SetFigFont{12}{14.4}{\rmdefault}{\mddefault}{\updefault}{\color[rgb]{0,0,0}$(1,1)$}%
}}}}
\end{picture}%

%% file: ladd2.pstex_t
\begin{picture}(0,0)%
\includegraphics{ladd2.pstex}%
\end{picture}%
\setlength{\unitlength}{3947sp}%
\begingroup\makeatletter\ifx\SetFigFont\undefined%
\gdef\SetFigFont#1#2#3#4#5{%
  \reset@font\fontsize{#1}{#2pt}%
  \fontfamily{#3}\fontseries{#4}\fontshape{#5}%
  \selectfont}%
\fi\endgroup%
\begin{picture}(3624,3624)(1189,-3973)
\put(2176,-1561){\makebox(0,0)[lb]{\smash{{\SetFigFont{12}{14.4}{\rmdefault}{\mddefault}{\updefault}{\color[rgb]{0,0,0}$\mathcal L^+$}%
}}}}
\end{picture}%

%% file: ladd3b.pstex_t
\begin{picture}(0,0)%
\includegraphics{ladd3b.pstex}%
\end{picture}%
\setlength{\unitlength}{3947sp}%
\begingroup\makeatletter\ifx\SetFigFont\undefined%
\gdef\SetFigFont#1#2#3#4#5{%
  \reset@font\fontsize{#1}{#2pt}%
  \fontfamily{#3}\fontseries{#4}\fontshape{#5}%
  \selectfont}%
\fi\endgroup%
\begin{picture}(3624,3918)(1189,-3973)
\put(2026,-1336){\makebox(0,0)[lb]{\smash{{\SetFigFont{12}{14.4}{\rmdefault}{\mddefault}{\updefault}{\color[rgb]{0,0,0}$\mathcal L^+_{(v,w)}$}%
}}}}
\put(3001,-211){\makebox(0,0)[lb]{\smash{{\SetFigFont{12}{14.4}{\rmdefault}{\mddefault}{\updefault}{\color[rgb]{0,0,0}$(v,w)$}%
}}}}
\end{picture}%

%% file: ladd4b.pstex_t
\begin{picture}(0,0)%
\includegraphics{ladd4b.pstex}%
\end{picture}%
\setlength{\unitlength}{3947sp}%
\begingroup\makeatletter\ifx\SetFigFont\undefined%
\gdef\SetFigFont#1#2#3#4#5{%
  \reset@font\fontsize{#1}{#2pt}%
  \fontfamily{#3}\fontseries{#4}\fontshape{#5}%
  \selectfont}%
\fi\endgroup%
\begin{picture}(3624,3657)(1189,-3973)
\put(2551,-811){\makebox(0,0)[lb]{\smash{{\SetFigFont{12}{14.4}{\rmdefault}{\mddefault}{\updefault}{\color[rgb]{0,0,0}$\mathcal L^+$}%
}}}}
\put(1951,-1861){\makebox(0,0)[lb]{\smash{{\SetFigFont{12}{14.4}{\rmdefault}{\mddefault}{\updefault}{\color[rgb]{0,0,0}$\mathcal H^+$}%
}}}}
\end{picture}%

%% file: ladd5.pstex_t
\begin{picture}(0,0)%
\includegraphics{ladd5.pstex}%
\end{picture}%
\setlength{\unitlength}{3947sp}%
\begingroup\makeatletter\ifx\SetFigFont\undefined%
\gdef\SetFigFont#1#2#3#4#5{%
  \reset@font\fontsize{#1}{#2pt}%
  \fontfamily{#3}\fontseries{#4}\fontshape{#5}%
  \selectfont}%
\fi\endgroup%
\begin{picture}(3624,3657)(1189,-3973)
\put(1876,-1486){\makebox(0,0)[lb]{\smash{{\SetFigFont{12}{14.4}{\rmdefault}{\mddefault}{\updefault}{\color[rgb]{0,0,0}$\mathcal L'^+$}%
}}}}
\put(2476,-511){\makebox(0,0)[lb]{\smash{{\SetFigFont{12}{14.4}{\rmdefault}{\mddefault}{\updefault}{\color[rgb]{0,0,0}$\mathcal J^+$}%
}}}}
\end{picture}%

%% file: arxives_v2.bbl
\begin{thebibliography}{11}

\bibitem{br81} W. Bruns, The Eisenbud-Evans generalized principal
  ideal theorem and determinantal ideals,
  Proc. Amer. Math. Soc. \textbf{83} (1981), no. 1, 19--24.

\bibitem{br88} W. Bruns, U. Vetter, Determinantal rings, Lecture Notes
  in Mathematics \textbf{1327}, Springer-Verlag, Berlin (1988).

\bibitem{br93} W. Bruns, J. Herzog, Cohen-Macaulay rings, Cambridge
  Studies in Adv. Math. \textbf{39}, Cambridge University Press,
  Cambridge (1993).

\bibitem{co94a} A. Conca, Symmetric ladders, Nagoya
  Math. J. \textbf{136} (1994), 35--56.

\bibitem{co94b} A. Conca, Divisor class group and canonical class of
  determinantal rings defined by ideals of minors of a symmetric
  matrix, Arch. Math. (Basel) \textbf{63} (1994), no. 3, 216--224.

\bibitem{co94c} A. Conca, Gr\"obner bases of ideals of minors of a
  symmetric matrix, J. Algebra \textbf{166} (1994), no. 2, 406--421.

\bibitem{de76} C. De Concini, C. Procesi, A characteristic free
  approach to invariant theory, Advances in Math. \textbf{21} (1976),
  no. 3, 330--354.

\bibitem{de82} C. De Concini, D. Eisenbud, C. Procesi, Hodge algebras,
  Ast\'erisque \textbf{91}, Soci\'et\'e Math\'ematique de France,
  Paris (1982). 

\bibitem{de08u} E. De Negri, E. Gorla, G-biliaison of ladder Pfaffian
  varieties, J. Algebra \textbf{321} (2009), no. 9, 2637--2649.

\bibitem{ea71} M. Hochster, J. Eagon, Cohen-Macaulay rings, invariant
  theory, and the generic perfection of determinantal loci,
  Amer. J. Math. \textbf{93} (1971), 1020--1058.

%

\bibitem{go07a} E. Gorla, The G-biliaison class of symmetric determinantal 
schemes, J. Algebra \textbf{310} (2007), no. 2, 880--902.

\bibitem{go07b} E. Gorla, Mixed ladder determinantal varieties
  from two-sided ladders, J. Pure Appl. Alg. \textbf{211} (2007),
  no. 2, 433--444.

\bibitem{go07u} E. Gorla, A generalized Gaeta's Theorem, Compositio
  Math. \textbf{144} (2008), no. 3, 689--704.

\bibitem{ha92} R. Hartshorne, Generalized divisors on Gorenstein
  curves and a theorem of Noether, J. Math. Kyoto Univ. \textbf{26}
  (1986), no. 3, 375--386.

\bibitem{ha94} R. Hartshorne, Generalized divisors on Gorenstein
  schemes, Proceedings of Conference on Algebraic Geometry and Ring
  Theory in honor of Michael Artin, Part III (Antwerp, 1992), K-Theory
  \textbf{8} (1994), no. 3, 287--339.

\bibitem{ha07} R. Hartshorne, Generalized divisors and biliaison,
  Illinois J. Math. \textbf{51} (2007), no. 1, 83--98.

\bibitem{he92} J. Herzog, N. V. Trung, Gr\"obner bases and
  multiplicity of determinantal and Pfaffian ideals,
  Adv. Math. \textbf{96} (1992),  no. 1, 1--37.

\bibitem{kl79} H. Kleppe, D. Laksov, The generic perfectness of
  determinantal schemes, Algebraic geometry (Proc. Summer Meeting,
  Univ. Copenhagen, Copenhagen, 1978),  244--252, Lecture Notes in
  Math. \textbf{732}, Springer, Berlin (1979).

\bibitem{kl01} J. O. Kleppe, J. C. Migliore, R. M. Mir\'o-Roig,
  U. Nagel, and C. Peterson, Gorenstein liaison, complete intersection
  liaison invariants and unobstructedness,
  Mem. Amer. Math. Soc. \textbf{154} (2001), no. 732.

\bibitem{mi98} J. Migliore, Introduction to Liaison Theory and
  Deficiency Modules, Birkh\"auser, Progress in Mathematics
  \textbf{165} (1998).

\bibitem{st90} B. Sturmfels, Gr\"obner bases and Stanley
  decompositions of determinantal rings, Math. Z. \textbf{205} (1990),
  no. 1, 137--144.

\end{thebibliography}
